\documentclass[a4,12pt,multicol,makeidx]{article}
\def\origin{
  \clearpage
\vskip-\baselineskip\vskip-\topskip%
  \vbox to 0pt{\vskip-1in%
    \hbox to 0pt{\hskip-1in%
      \hbox to 0pt{\vrule width 1cm height .4pt depth 0mm\hss}%
      \vbox to 0pt{\hrule width .4pt height 0pt depth 1cm\vss}%
    \hss}%
  \vss}
  \vskip-\baselineskip
  \vbox to 0pt{\vskip-1in\vskip3cm%
    \hbox to 0pt{\hskip-1in\hskip3cm%
      \hbox to 0pt{\hss\vrule width 2cm height .4pt depth 0mm\hss}%
      \vbox to 0pt{\vss\hrule width .4pt height 1cm depth 1cm\vss}%
    \hss}%
  \vss}%
\vskip5mm\hskip10mm (3cm,3cm)
}%
              
    \def\l{\lambda}

\def\n{\nabla}

\def\n{\nabla} 
\def\p{\partial}  

\newenvironment{theorem}{%
\par \bigskip \it}{%
\bigskip \par}

\makeindex
\title{An Interior Gradient Estimate for a class of Second Order Partial Differential Inequalities. 
}
\author{Mariko Arisawa\\ Department of Mathematics
\\University of  Wisconsin-Madison\\
Madison, WI 53706, U.S.A.\\
}
\date{}
\begin{document}
\maketitle
\bigskip

\section{Introduction.}	

$\quad$ In this paper, we show a uniform gradient estimate for functions $u$ defined on a domain $\Omega$ in Euclidian space ${\bf R^N}$ which satisfy a system of $N$ second order partial differential inequalities of the following form
\begin{equation}
	-\frac{\p}{\p x_i}(\sum_{j=1}^{N}A_{ij}\frac{\p u}{\p x_j})(x)+ \sum_{j=1}^{N}b_{ij}\frac{\p u}{\p x_j}(x)\leq C
	\quad x\in  \Omega,\quad 1\leq i\leq N. 
\end{equation}
The type of this system we concern is discussed in our main results Theorems 1.1, 1.2, 1.3, and 1.4 by some conditions for the coefficient matrix $(A_{ij})_{1\leq i,j\leq N}$ of second order terms and for the coefficient matrix $(b_{ij})_{1\leq i,j\leq N}$ of first order terms. \\
$\quad$ Our method relies essentially on the structure of the system (1), and we shall prove an interior uniform gradient estimate for $u\in C^2(\Omega)$ under the assumption that $u$ is uniformly bounded. It is worth remarking here that the technical difficulty increases with the dimension $N$, and we only give the result for the cases of $N=2,3$. (Theorems 1.1, 1.2.) On the other hand, if we assume that $u$ has a compact support in $\Omega$, or that $u$ is periodic, the boundedness of $u$ is not necessary and the same result holds with a simpler assumption on $(A_{ij})_{1\leq i,j\leq N}$ for any dimensions. (Theorems 1.3, 1.4.) 
We may now state our main theorems. \\

\leftline{\bf Theorem 1.1.}
\begin{theorem}
Let $\Omega$ be a domain in ${\bf R^2}$, let $A=$$(A_{ij})_{1\leq i,j\leq 2}$, where $A_{ij}$$=A_{ij}(x_1,x_2)$$\in W^{1,\infty}(\Omega)$ are real valued functions defined in $(x_1,x_2)\in \Omega$ which satisfy the following conditions. 
\begin{equation}
	\sup_{(x_1,x_2)\in \Omega} |A_{ij}(x_1,x_2)|\leq C_1,\quad \sup_{(x_1,x_2)\in \Omega} |\n A_{ij}(x_1,x_2)|\leq C_1,
	\quad 1\leq i,j\leq 2,
\end{equation}
\begin{equation}
	|det A|^{-1} = |A_{11}A_{22}-A_{12}A_{21}|^{-1} \leq C_2,
\end{equation}
\begin{equation}
	A_{11},\quad A_{22}\neq 0,
\end{equation}
where $C_1$, $C_2>0$ are constants. Let $b_{ij}$$=b_{ij}(x_1,x_2)$$\in W^{1,\infty}(\Omega)$, ($1\leq i,j\leq 2$) be real valued functions defined in $(x_1,x_2)$$\in \Omega$ which satisfy the following conditions.  
\begin{equation}
	\sup_{(x_1,x_2)\in \Omega} |b_{ij}(x_1,x_2)|\leq C_3,\quad \sup_{(x_1,x_2)\in \Omega} |\n b_{ij}(x_1,x_2)|\leq C_3,
	\quad 1\leq i,j\leq 2,
\end{equation}
where $C_3>0$ is a constant. Suppose that a real valued function $u(x_1,x_2)$$\in C^2(\Omega)$ satisfies the following inequalities 
\begin{equation}
	-\sum_{j=1}^{2}\{ \frac{\p}{\p x_i}(A_{ij}\frac{\p u}{\p x_j})+ b_{ij}\frac{\p u}{\p x_j} \}(x_1,x_2)\leq C_4
	\quad \hbox{in}\quad (x_1,x_2)\in  \Omega,\quad i=1,2. 
\end{equation}
\begin{equation}
	\sup_{(x_1,x_2)\in \Omega} |u|\leq C_5,
\end{equation}
where $C_4$, $C_5>0$ are constants. Then, for any $(y_1,y_2)$$\in \Omega$, any $\delta>0$ such that  
$$
	K'=[y_1-2\delta,y_1+2\delta]\times [y_2-2\delta,y_2+2\delta] \subset \Omega,
$$
for $K=[y_1-\delta,y_1+\delta]\times [y_2-\delta,y_2+\delta]$, there exists a constant $\overline{C}>0$ depending on $\delta$, matrices $(A_{ij})$, $(b_{ij})$, and constants $C_4$, $C_5$ such that 
\begin{equation}
	\sup_{(x_1,x_2)\in K} |\n u|\leq \overline{C},
\end{equation}
$$
	\overline{C}= 
	O(\frac{1}{\delta}). 
$$
\end{theorem}

\leftline{\bf Theorem 1.2.}
\begin{theorem}
Let $\Omega$ be a domain in ${\bf R^3}$, let $A$$=(A_{ij})_{1\leq i,j\leq 3}$, where $A_{ij}$ are constants which satisfy the following conditions. 
\begin{equation}
	\sup_{(x_1,x_2,x_3)\in \Omega} |A_{ij}(x_1,x_2,x_3)|\leq C_1
	\quad 1\leq i,j\leq 3,
\end{equation}
\begin{equation}
	|det A|^{-1}  \leq C_2,
\end{equation}
\begin{equation}
	(A_{11}A_{22}-A_{12}A_{21})(A_{11}A_{33}-A_{13}A_{31})(A_{22}A_{33}-A_{23}A_{32}) \neq 0,
\end{equation}
\begin{equation}
	A_{11}A_{22}A_{33}\neq 0,
\end{equation}
\begin{equation}
	(A_{11}A_{22}+A_{21}A_{12})(A_{11}A_{22}-3A_{21}A_{12}) > 0,
\end{equation}
\begin{equation}
	(A_{11}A_{33}+A_{31}A_{13})(A_{11}A_{33}-3A_{31}A_{13}) > 0,
\end{equation}
\begin{equation}
	(A_{22}A_{33}+A_{23}A_{32})(A_{22}A_{33}-3A_{23}A_{32}) > 0,
\end{equation}
where $C_1$, $C_2>0$ are constants, let $b_{ij}$ ($1\leq i,j\leq 3$) be constants which satisfy 
\begin{equation}
	\sup_{(x_1,x_2,x_3)\in \Omega} |b_{ij}|\leq C_3
	\quad 1\leq i,j\leq 3,
\end{equation}
where $C_3>0$ is a constant. Suppose that a real valued function $u(x_1,x_2,x_3)$$\in C^2(\Omega)$ satisfies the following inequalities 
\begin{equation}
	-\sum_{j=1}^{3} \{ \frac{\p}{\p x_i}(A_{ij}\frac{\p u}{\p x_j})+ b_{ij}\frac{\p u}{\p x_j}\} (x_1,x_2,x_3)\leq C_4, 
\end{equation}
$$
	\hbox{in}\quad  (x_1,x_2,x_3) \in  \Omega,\quad 1\leq \forall i\leq 3, 
$$
\begin{equation}
	\sup_{(x_1,x_2,x_3)\in \Omega} |u|\leq C_5,
\end{equation}
where $C_4$, $C_5$ are constants. Then, for any $(y_1,y_2,y_3)$$\in \Omega$ and for any $\delta>0$ such that 
$$
	\delta |\frac{A_{ij}b_{ji}}{A_{ii}A_{jj}-A_{ij}A_{ji}+2|A_{ij}A_{ji}|}| (1+|\frac{A_{ij}A_{ji}}{A_{ii}A_{jj}}|) <\frac{1}{8} \quad 
	1\leq i\neq j\leq 3, 
$$
\begin{equation}\quad
\end{equation}
$$
	\delta |\frac{A_{ij}b_{ij}}{A_{ii}A_{jj}-A_{ij}A_{ji}-2|A_{ij}A_{ji}|} |(1+|\frac{A_{ij}A_{ji}}{A_{ii}A_{jj}}|) <\frac{1}{8} \quad 
	1\leq i\neq j\leq 3, 
$$
$$
	K'=[y_1-2\delta,y_1+2\delta]\times [y_2-2\delta,y_2+2\delta]\times [y_3-2\delta,y_3+2\delta] \subset \Omega,
$$
for $K=[y_1-\delta,y_1+\delta]\times [y_2-\delta,y_2+\delta]\times [y_3-\delta,y_3+\delta]$, there exists a constant $\overline{C}>0$ depending on $\delta$, 
matrices $(A_{ij})$, $(b_{ij})$, and constants $C_4$, $C_5$ such that 
\begin{equation}
	\sup_{(x_1,x_2)\in K} |\n u|\leq \overline{C},
\end{equation}
$$
	\overline{C}= O(\frac{1}{\delta}). 
$$
\end{theorem}

\leftline{\bf Theorem 1.3.}
\begin{theorem}
Let $\Omega$ be a domain in ${\bf R^N}$, let $A$$=(A_{ij})_{1\leq i,j\leq N}$, where $A_{ij}$$=A_{ij}$$\in L^{\infty}(\Omega)$ ($1\leq i,j\leq N$) real valued functions defined in $x\in \Omega$ which satisfy the following conditionzs. 
\begin{equation}
	\sup_{x \in \Omega} |A_{ij}(x)|\leq C_1
	\quad 1\leq i,j\leq N,
\end{equation}
\begin{equation}
	|det A|^{-1}  \leq C_2,
\end{equation}
where $C_1$, $C_2$ are constants. Suppose that a real valued function $u(x)\in C^2(\Omega)$ such that $\hbox{supp} u\subset \subset \Omega$ satisfies the following inequalities
\begin{equation}
	-\frac{\p}{\p x_i}(\sum_{j=1}^{N}A_{ij}\frac{\p u}{\p x_j})(x)\leq C_3 \quad \hbox{in} 
	\quad x\in  \Omega,\quad 1\leq i\leq N, 
\end{equation}
where $C_3>0$ is a constant. Then, there exists a constant $\overline{C}>0$ depending on the matrix $(A_{ij})$ and the constant $C_3>0$ such that 
\begin{equation}
	\sup_{x\in \Omega} |\n u(x)|\leq \overline{C}. 
\end{equation}
\end{theorem}

\leftline{\bf Theorem 1.4.}
\begin{theorem}
 Let $\Omega$ be an $N$ dimensional torus ${\bf T^N}$$={\bf R^N/Z^N}$$=[0,1]^N$, let $A$$=(A_{ij})$, where $A_{ij}$$=A_{ij}(x)$$\in L^{\infty}(\Omega)$ ($1\leq i,j\leq N$) real valued periodic functions defined in $x\in \Omega$ which satisfy the following conditions. 
\begin{equation}
	\sup_{x \in \Omega} |A_{ij}(x)|\leq C_1
	\quad 1\leq i,j\leq N,
\end{equation}
\begin{equation}
	|det A|^{-1}  \leq C_2,
\end{equation}
where $C_1$, $C_2>0$ are constants. Suppose that a real valued function $u(x)\in C^2(\Omega)$ is periodic in $\Omega$ and satisfies the following inequalities 
\begin{equation}
	-\frac{\p}{\p x_i}(\sum_{j=1}^{N}A_{ij}\frac{\p u}{\p x_j})(x)\leq C_3 \quad \hbox{in} 
	\quad x\in  \Omega,\quad 1\leq i\leq N, 
\end{equation}
where $C_3>0$ is a constant. Then, there exists a constant $\overline{C}>0$ depending on the matrix $(A_{ij})$ and the constant $C_3>0$ such that  
\begin{equation}
	\sup_{x\in \Omega} |\n u(x)|\leq \overline{C}. 
\end{equation}
\end{theorem}

$\quad$ If we do not assume either the condition $\hbox{supp}u$$\subset \subset \Omega$ in Theorem 1.3, or the periodicity in Theorem 1.4, we need more restrictive conditions for the matrix $A$ as in Theorems 1.1, 1.2. The following counter examples show the contrast between the case of Theorems 1.1, 1.2 and the case of Theorems 1.3, 1.4. \\

\leftline{\bf Example 1.5.}
 Let $N=2$, and let $A=$$(A_{ij})_{1\leq i,j\leq 2}$ be the matrix with $A_{11}$$=A_{22}$$=0$, $A_{21}$$=A_{12}=1$. ($\hbox{det} A\neq 0$.) Consider any functions $u(x_1,x_2)$$\in C^2(\Omega)$ which satisfy the following inequality in $(x_1,x_2)\in \Omega$. 
$$
	-\frac{\p^2 u}{\p x_1 \p x_2}\leq C_0. 
$$
Then, if $\hbox{supp}u$$\subset\subset$$\Omega$, from Theorem 1.3 $|\n u|\leq C$, where the constant $C>0$ depends only on the matrix $A$ and $C_0$. 
However, if we take the function $u(x_1,x_2)=\psi(x_1)$ with arbitrary $\psi\in C^2({\bf R})$ such that $\hbox{supp}u\cap \p\Omega$$\neq \emptyset$, 
although $u$ satisfies the above partial differential inequality, $|\n u|$ is not bounded in general. \\

\leftline{\bf Example 1.6.}
 Let $N=3$, and let $A=$$(A_{ij})_{1\leq i,j\leq 3}$ be the matrix with $A_{11}$ $=A_{12}$ $=A_{21}$ $=A_{23}$ $=0$, $A_{13}$ $=A_{22}$ 

$=A_{31}=$$A_{32}$ $=A_{33}$ $=1$. 

($\hbox{det} A \neq 0$) 
Consider any functions $u(x_1,x_2,x_3)$$\in C^2(\Omega)$ which satisfy the following inequality in $(x_1,x_2,x_3)\in \Omega$. 

$$
	-\frac{\p^2 u}{\p x_1 \p x_3}\leq C_0. 
$$
$$
	-\frac{\p^2 u}{\p x_2^2}\leq C_0. 
$$
$$
	-\frac{\p }{\p x_3}(\frac{\p u}{\p x_1} + \frac{\p u}{\p x_2} + \frac{\p u}{\p x_3}) \leq C_0. 
$$
Then, if $\hbox{supp}u$$\subset\subset$$\Omega$, from Theorem 1.3 $|\n u|\leq C$, where the constant $C>0$ depends only on the matrix $A$ and $C_0$.
However, if we take the function $u(x_1,x_2,x_3)=\psi(x_1,x_2)$ with arbitrary $\psi\in C^2({\bf R^2})$ such that $-\frac{\p^2 \psi}{\p x_2^2}<0$, $\hbox{supp}u\cap \p\Omega$$\neq \emptyset$, 
although $u$ satisfies the above partial differential inequalities, $|\n u|$ is not bounded in general. \\

$\quad$ Finally, we shall give an example of second-order degenerate elliptic partial differential equation whose regularity of the solution can be shown by our results. \\

\leftline{\bf Example 1.7.}
 Let $N=2,3$, $\Omega$ be a bounded domain in ${\bf R^N}$, $\l>0$, and suppose that $u_{\l}$ is a solution of the following problem
\begin{equation}
	\l u_{\l}(x) + \sup_{1\leq i\leq N} \{ -\frac{\p^2 u_{\l}}{\p x_i^2}(x) \} -V(x)=0 \quad x\in \Omega,  
\end{equation}
with either Neumann B.C., or State constraint B.C., where $V(x)$ is a Lipschitz continuous function defined in $\Omega$. Then, for any interior domain $\Omega_0$$\subset\subset$$\Omega$, $\l u_{\l}(x)$ is Lipschitz continuous in $x\in \Omega_0$ uniformly with respect to $\l>0$. If we assume that a solution $u_{\l}$ of (29) satisfies Periodic B.C., then $u_{\l}(x)$ is Lipschitz continuous in $x\in \Omega$ uniformly with respect to $\l$. \\

$\quad$ The proof of Example 1.7 will be given later in this paper. The partial differential equation (29) corresponds to an optimal control problem for a stochastic system which equips $N$ controls of one dimensional diffusions at each state. (Remark that $\sum_{i=1}^{N}(-\frac{\p^2}{\p x_i^2})=-\Delta$.) In view of the degeneracy of (29), we cannot apply the usual regularity theory for uniformly elliptic operators (Gilbarg-Trudinger [5], Caffarelli-Cabre [2]) to study the regularity of the solution $u_{\l}$.\\

$\quad$ The plan of this paper is as follows. Theorem 1.1 is proved in \S 2; Theorem 1.2 is proved in \S 3; and the proofs of Theorems 1.3, 1.4 and Example 1.7 are given in \S 4. Throughout in the present paper, we conserve the letter $C>0 $  to denote the constants which depend on constants $C_i$ and matrices $(A_{ij})$, $(b_{ij})$ in the Theorems 1.1-1.4. \\

\section{Proof of Theorem 1.1.}

\leftline{\bf Lemma 2.1.}
\begin{theorem}
 Let $\phi$ be an arbitrarily fixed real valued twice differentiable function defined on the interval $[-2\delta,2\delta]$ such that $0\leq \phi\leq 1$, $\hbox{supp} \phi$$\subset\subset$$(-2\delta,2\delta)$, $\phi$ is even and 
\begin{equation}
	\phi=1\quad \hbox{on}\quad [-\delta,\delta],\quad \phi'\geq 0\quad \hbox{on}\quad [-2\delta,0]. 
\end{equation}
Then, the function $u$ in Theorem 1.1 satisfies the following inequalities : for any $(\hat{x}_1, x_2)\in K$, 
\begin{equation}
	(A_{11}\frac{\p u}{\p x_1})(\hat{x}_1, x_2)
\end{equation}
$$
	\leq (resp. \geq) \quad -(A_{12}\frac{\p u}{\p x_2})(\hat{x_1},x_2) + \frac{1}{2\delta} \int_{-\delta}^{\delta} (A_{12}\frac{\p u}{\p x_2})({x'_1},x_2) dx'_1
$$
$$
	-\frac{1}{2\delta} \int_{-\delta}^{\delta}  \int_{\hat{x}_1}^{{x}'_1}    (b_{12}\frac{\p u}{\p x_2})({x''_1},x_2) dx''_1dx'_1
$$
$$
	-(resp. +) \int_{-2\delta}^{2\delta} \phi'(x'_1)(A_{12}\frac{\p u}{\p x_2})({x'_1},x_2) dx'_1
$$
$$
	-(resp. +) \int_{-2\delta}^{2\delta} \phi(x'_1)(b_{12}\frac{\p u}{\p x_2})({x'_1},x_2) dx'_1+(resp.-)\quad C; 
$$
for any $({x_1},\hat{x}_2)\in K$, 
\begin{equation}
	(A_{22}\frac{\p u}{\p x_2})(x_1,\hat{x}_2)
\end{equation}
$$
	\leq (resp. \geq) \quad -(A_{21}\frac{\p u}{\p x_1})({x_1},\hat{x}_2) + \frac{1}{2\delta} \int_{-\delta}^{\delta} (A_{21}\frac{\p u}{\p x_1})({x_1},x'_2) dx'_2
$$
$$
	-\frac{1}{2\delta} \int_{-\delta}^{\delta}  \int_{\hat{x}_2}^{{x}'_2}    (b_{21}\frac{\p u}{\p x_1})({x_1},x''_2) dx''_2dx'_2
$$
$$
	-(resp. +) \int_{-2\delta}^{2\delta} \phi'(x'_2)(A_{21}\frac{\p u}{\p x_1})({x_1},x'_2) dx'_2
$$
$$
	-(resp. +) \int_{-2\delta}^{2\delta} \phi(x'_2)(b_{21}\frac{\p u}{\p x_1})({x_1},x'_2) dx'_2+(resp.-)\quad C; 
$$
\end{theorem}

\leftline{\bf Lemma 2.2.}
\begin{theorem}
 For the terms in (31) the following estimates hold 
$$
	|\int_{-\delta}^{\delta}  \int_{\hat{x}_1}^{{x}'_1}    (b_{12}\frac{\p u}{\p x_2})({x''_1},x_2) dx''_1dx'_1|\leq C,
$$
$$
	|\int_{-2\delta}^{2\delta} \phi'(x'_1)(A_{12}\frac{\p u}{\p x_2})({x'_1},x_2) dx'_1|\leq C,
$$
\begin{equation}
\quad
\end{equation}
$$
	|\int_{-2\delta}^{2\delta} \phi(x'_1)(b_{12}\frac{\p u}{\p x_2})({x'_1},x_2) dx'_1|\leq C,
$$
$$
	|\int_{-\delta}^{\delta}  (A_{12}\frac{\p u}{\p x_2})({x'_1},x_2) dx'_1|\leq C. 
$$
For the terms in (32), the following estimates hold
$$
	|\int_{-\delta}^{\delta}  \int_{\hat{x}_2}^{{x}'_2}    (b_{21}\frac{\p u}{\p x_1})({x_1},x''_2) dx''_2dx'_2|\leq C,
$$
$$
	|\int_{-2\delta}^{2\delta} \phi'(x'_2)(A_{21}\frac{\p u}{\p x_1})({x_1},x'_2) dx'_2|\leq C,
$$
\begin{equation}
\quad
\end{equation}
$$
	|\int_{-2\delta}^{2\delta} \phi(x'_2)(b_{21}\frac{\p u}{\p x_1})({x_1},x'_2) dx'_2|\leq C,
$$
$$
	|\int_{-\delta}^{\delta}  (A_{21}\frac{\p u}{\p x_1})({x_1},x'_2) dx'_2|\leq C. 
$$
\end{theorem}

$\quad$ We temporarily admit Lemmas 2.1, 2.2, and give the proof of Theorem 1.1. Let us remark that (31)-(33), (32)-(34) lead the following estimates for any $(\hat{x}_1,x_2)$, $({x}_1,\hat{x}_2)$$\in K$.  
\begin{equation}
	-(A_{12}\frac{\p u}{\p x_2})(\hat{x}_1,x_2) - C\leq (A_{11}\frac{\p u}{\p x_1})(\hat{x}_1,x_2)
	\leq -(A_{12}\frac{\p u}{\p x_2})(\hat{x}_1,x_2)+C 
\end{equation}
\begin{equation}
	-(A_{21}\frac{\p u}{\p x_1})({x}_1,\hat{x}_2) - C\leq (A_{22}\frac{\p u}{\p x_2})({x}_1,\hat{x}_2)
	\leq -(A_{21}\frac{\p u}{\p x_1})({x}_1,\hat{x}_2)+C 
\end{equation}
From (4), (36), 
$$
	-(A_{12}\frac{\p u}{\p x_2})(\hat{x}_1,x_2) = -(\frac{A_{12}}{A_{22}} A_{22} \frac{\p u}{\p x_2})(\hat{x}_1,x_2) 
$$
$$
	\leq (resp.\quad\geq )\quad (\frac{A_{21}A_{12}}{A_{22}}  \frac{\p u}{\p x_1})(\hat{x}_1,x_2) + (resp. -)C, 
$$
and by inserting this into (35) we get for any $ (\hat{x}_1,x_2)$$\in K$, 
$$
	|(\frac{A_{11}A_{22}-A_{21}A_{12}}{A_{22}}  \frac{\p u}{\p x_1})(\hat{x}_1,x_2) |\leq C, 
$$
where $C>0$ is a constant. From the assumption (3), (4), we have the bound for $\frac{\p u}{\p x_1}$, and the same discussion leads the bound for $\frac{\p u}{\p x_2}$. Therefore, (8) is proved. \\

$\quad$ Now, we shall give the proof of the Lemmas 2.1, 2.2. \\

{\bf Proof of Lemma 2.1.}\\
 $\quad$ We only give the proof of (31); (32) will be obtained by the same way. First, from (6) for any $(x_1,x_2)$, $(\hat{x}_1,x_2)$$\in K$, since 
$$
	-(A_{11}\frac{\p u}{\p x_1}+ A_{12}\frac{\p u}{\p x_2})({x}_1,x_2)+ \int_{\hat{x}_1}^{x_1} (b_{12}\frac{\p u}{\p x_2})(x'_1,x_2) dx'_1
$$
$$
	= -(A_{11}\frac{\p u}{\p x_1}+ A_{12}\frac{\p u}{\p x_2})(\hat{x}_1,x_2)
$$
$$
	\int_{\hat{x}_1}^{x_1} \{-\frac{\p}{\p x'_1} (A_{11}\frac{\p u}{\p x'_1}+ A_{12}\frac{\p u}{\p x_2})({x}'_1,x_2)
    + (b_{12}\frac{\p u}{\p x_2})(x'_1,x_2) \}dx'_1
$$
by using (7), (2), (5), the following holds with $\phi$ stated in Lemma 2.1. 
$$
	|-(A_{11}\frac{\p u}{\p x_1}+ A_{12}\frac{\p u}{\p x_2})({x}_1,x_2)+ \int_{\hat{x}_1}^{x_1} (b_{12}\frac{\p u}{\p x_2})(x'_1,x_2) dx'_1
$$
$$
	+ (A_{11}\frac{\p u}{\p x_1}+ A_{12}\frac{\p u}{\p x_2})(\hat{x}_1,x_2)|
$$
$$
	\leq \hbox{sgn}(x_1-\hat{x}_1) \int_{\hat{x}_1}^{x_1} | - \frac{\p}{\p x'_1} (A_{11}\frac{\p u}{\p x'_1}+ A_{12}\frac{\p u}{\p x_2})({x}'_1,x_2)
$$
$$
    + (b_{11}\frac{\p u}{\p x'_1}+ b_{12}\frac{\p u}{\p x_2})(x'_1,x_2) |dx'_1 +C 
$$
$$
	= \hbox{sgn}(x_1-\hat{x}_1)   \int_{\hat{x}_1}^{x_1} | C_0+ \frac{\p}{\p x'_1} (A_{11}\frac{\p u}{\p x'_1}+ A_{12}\frac{\p u}{\p x_2})({x}'_1,x_2)
$$
$$
	-(b_{11}\frac{\p u}{\p x'_1}+ b_{12}\frac{\p u}{\p x_2})(x'_1,x_2) -C_0|dx'_1 +C 
$$
$$
	\leq \hbox{sgn}(x_1-\hat{x}_1) \times 
$$
$$
	 \int_{\hat{x}_1}^{x_1}  
		C_0+ \frac{\p}{\p x'_1} (A_{11}\frac{\p u}{\p x'_1}+ A_{12}\frac{\p u}{\p x_2})({x}'_1,x_2)
		-(b_{11}\frac{\p u}{\p x'_1}+ b_{12}\frac{\p u}{\p x_2})(x'_1,x_2) dx'_1 +C 
$$
$$
	=\hbox{sgn}(x_1-\hat{x}_1) 
	\int_{\hat{x}_1}^{x_1}  \phi(x'_1)\{
		C_0+ \frac{\p}{\p x'_1} (A_{11}\frac{\p u}{\p x'_1}+ A_{12}\frac{\p u}{\p x_2})({x}'_1,x_2)
$$
$$
	-(b_{11}\frac{\p u}{\p x'_1}+ b_{12}\frac{\p u}{\p x_2})(x'_1,x_2) \}dx'_1 +C 
$$
$$
	\leq \int_{-2\delta}^{2\delta}   \phi(x'_1)\{ \frac{\p}{\p x'_1} (A_{11}\frac{\p u}{\p x'_1}+ A_{12}\frac{\p u}{\p x_2})({x}'_1,x_2)
	-(b_{11}\frac{\p u}{\p x'_1}+ b_{12}\frac{\p u}{\p x_2})(x'_1,x_2) \}dx'_1 +C 
$$
$$
	\leq - \int_{-2\delta}^{2\delta}   \phi(x'_1) (A_{11}\frac{\p u}{\p x'_1}+ A_{12}\frac{\p u}{\p x_2})({x}'_1,x_2)dx'_1
	-\int_{-2\delta}^{2\delta}  \phi(x'_1)(b_{12}\frac{\p u}{\p x_2})(x'_1,x_2) dx'_1 +C 
$$
$$
	= - \int_{-2\delta}^{2\delta}   \phi'(x'_1)(A_{12}\frac{\p u}{\p x_2})({x}'_1,x_2)dx'_1
	-\int_{-2\delta}^{2\delta}  \phi(x'_1)(b_{12}\frac{\p u}{\p x_2})(x'_1,x_2) dx'_1 +C. 
$$
Hence, for any $(x_1,x_2)$, $(\hat{x}_1,x_2)$$\in K$ we get the following inequalities. 
\begin{equation}
	(A_{11}\frac{\p u}{\p x_1}+ A_{12}\frac{\p u}{\p x_2})({x}_1,x_2)
\end{equation}
$$
	\leq (resp.\geq) \quad (A_{11}\frac{\p u}{\p x_1}+ A_{12}\frac{\p u}{\p x_2})(\hat{x}_1,x_2)
	+\int_{\hat{x}_1}^{x_1} b_{12}\frac{\p u}{\p x_2}({x}'_1,x_2)dx'_1
$$
$$
	-(resp.+) \int_{-2\delta}^{2\delta}  \phi'(x'_1)(A_{12}\frac{\p u}{\p x_2})(x'_1,x_2) dx'_1
$$
$$
	-(resp.+) \int_{-2\delta}^{2\delta}  \phi(x'_1)(b_{12}\frac{\p u}{\p x_2})(x'_1,x_2) dx'_1 +(resp.-)\quad C. 
$$
Next, we integrate both hand sides of the above inequalities with respect to $x_1$ on $[-\delta,\delta]$, then devide the obtained result by $2\delta$ and we have the following. 
\begin{equation}
	\frac{1}{2\delta} \int_{-\delta}^{\delta} (A_{12}\frac{\p u}{\p x_2})({x}'_1,x_2)dx'_1
\end{equation}
$$
	\leq \quad (resp.\geq)\quad (A_{11}\frac{\p u}{\p x_1}+ A_{12}\frac{\p u}{\p x_2})(\hat{x}_1,x_2)
$$
$$
	+ \frac{1}{2\delta} \int_{-\delta}^{\delta} \int_{\hat{x}_1}^{x'_1} (b_{12}\frac{\p u}{\p x_2})({x}''_1,x_2)dx''_1dx'_1
$$
$$
	-(resp.+) \int_{-2\delta}^{2\delta}  \phi'(x'_1)(A_{12}\frac{\p u}{\p x_2})(x'_1,x_2) dx'_1
$$
$$
	-(resp.+) \int_{-2\delta}^{2\delta}  \phi(x'_1)(b_{12}\frac{\p u}{\p x_2})(x'_1,x_2) dx'_1 +(resp.-)\quad C. 
$$
From (38), (31) holds clearly.\\

{\bf Proof of Lemma 2.2.}\\
$\quad$ In the two dimensional case, the argument is easy. In fact, to show the first inequalities in (33), we multiply both hand sides of the inequalities (32) by 
$(A_{22}^{-1}b_{12})(x_1,\hat{x}_2)$ and then integrate the result first with respect to $x_1$ on $[\hat{x}_1,x'_1]$, then with respect to $x'_1$ on $[-\delta,\delta]$, 
which leads the conclusion because of (2), (4), (5), (7). The other estimates can be obtained similarly, which we do not repeat here. \\

\section{Proof of Theorem 1.2.}
\leftline{\bf Lemma 3.1.}
\begin{theorem}
$\quad$ Let $\phi$ be an arbotrarily fixed real valued twice differentiable function  defined on the interval $[-2\delta,2\delta]$ such that $0\leq \phi\leq 1$, 
$\hbox{supp}\phi$$\subset\subset$$(-2\delta,2\delta)$, $\phi$ is even and 
\begin{equation}
	\phi=1\quad \hbox{on}\quad [-\delta,\delta], \quad \phi'\geq 0\quad \hbox{on}\quad [-2\delta,0]. 
\end{equation}
Then, the function $u$ in Theorem 1.2 satisfies the following inequalities: for any $(\hat{x}_1,x_2,x_3)$$\in K$, 
\begin{equation}
	(A_{11}\frac{\p u}{\p x_1})(\hat{x}_1,x_2,x_3)
\end{equation}
$$
	\leq (resp.\geq)\quad - (A_{12}\frac{\p u}{\p x_2}+ A_{13}\frac{\p u}{\p x_3})(\hat{x}_1,x_2,x_3)
$$
$$
	+ \frac{1}{2\delta} \int_{-\delta}^{\delta} (A_{12}\frac{\p u}{\p x_2}+A_{13}\frac{\p u}{\p x_3})({x}'_1,x_2,x_3)dx'_1
$$
$$
	- \frac{1}{2\delta} \int_{-\delta}^{\delta} \int_{\hat{x}_1}^{x'_1} (b_{12}\frac{\p u}{\p x_2} +b_{13}\frac{\p u}{\p x_3} )({x}''_1,x_2,x_3)dx''_1dx'_1
$$
$$
	-(resp.+) \int_{-2\delta}^{2\delta}  \phi'(x'_1)(A_{12}\frac{\p u}{\p x_2}+ A_{13}\frac{\p u}{\p x_3})(x'_1,x_2,x_3) dx'_1
$$
$$
	-(resp.+) \int_{-2\delta}^{2\delta}  \phi(x'_1)(b_{12}\frac{\p u}{\p x_2}+ b_{13}\frac{\p u}{\p x_3})(x'_1,x_2,x_3) dx'_1 +(resp.-)\quad C, 
$$
for any $({x}_1,\hat{x}_2,x_3)\in K$, 
\begin{equation}
	(A_{22}\frac{\p u}{\p x_2})({x}_1,\hat{x}_2,x_3)
\end{equation}
$$
	\leq (resp.\geq)\quad - (A_{21}\frac{\p u}{\p x_1}+ A_{23}\frac{\p u}{\p x_3})({x}_1,\hat{x}_2,x_3)
$$
$$
	+ \frac{1}{2\delta} \int_{-\delta}^{\delta} (A_{21}\frac{\p u}{\p x_1}+A_{23}\frac{\p u}{\p x_3})({x}_1,x'_2,x_3)dx'_2
$$
$$
	- \frac{1}{2\delta} \int_{-\delta}^{\delta} \int_{\hat{x}_2}^{x'_2} (b_{21}\frac{\p u}{\p x_1} +b_{23}\frac{\p u}{\p x_3} )(x_1,x''_2,x_3)dx''_2dx'_2
$$
$$
	-(resp.+) \int_{-2\delta}^{2\delta}  \phi'(x'_2)(A_{21}\frac{\p u}{\p x_1}+ A_{23}\frac{\p u}{\p x_3})(x_1,x'_2,x_3) dx'_2
$$
$$
	-(resp.+) \int_{-2\delta}^{2\delta}  \phi(x'_2)(b_{21}\frac{\p u}{\p x_1}+ b_{23}\frac{\p u}{\p x_3})(x_1,x'_2,x_3) dx'_2 +(resp.-)\quad C. 
$$
for any $({x}_1,x_2,\hat{x}_3)$$\in K$, 
\begin{equation}
	(A_{33}\frac{\p u}{\p x_3})({x}_1,{x}_2,\hat{x}_3)
\end{equation}
$$
	\leq (resp.\geq)\quad - (A_{31}\frac{\p u}{\p x_1}+ A_{32}\frac{\p u}{\p x_2})({x}_1,{x}_2,\hat{x}_3)
$$
$$
	+ \frac{1}{2\delta} \int_{-\delta}^{\delta} (A_{31}\frac{\p u}{\p x_1}+A_{32}\frac{\p u}{\p x_2})({x}_1,x_2,x'_3)dx'_3
$$
$$
	- \frac{1}{2\delta} \int_{-\delta}^{\delta} \int_{\hat{x}_3}^{x'_3} (b_{31}\frac{\p u}{\p x_1} +b_{32}\frac{\p u}{\p x_2} )(x_1,x_2,x''_3)dx''_3dx'_3
$$
$$
	-(resp.+) \int_{-2\delta}^{2\delta}  \phi'(x'_3)(A_{31}\frac{\p u}{\p x_1}+ A_{32}\frac{\p u}{\p x_2})(x_1,x_2,x'_3) dx'_3
$$
$$
	-(resp.+) \int_{-2\delta}^{2\delta}  \phi(x'_3)(b_{31}\frac{\p u}{\p x_1}+ b_{32}\frac{\p u}{\p x_2})(x_1,x_2,x'_3) dx'_3 +(resp.-)\quad C. 
$$
\end{theorem}

\leftline{\bf Lemma 3.2.}
\begin{theorem}
$\quad$ Let us denote $({x}'_1,x_2,x_3)$$=(x'_1)$, $({x}_1,x'_2,x_3)$$=(x'_2)$, $({x}_1,x_2,x'_3)$$=(x'_3)$, and 
$({x}''_1,x_2,x_3)$$=(x''_1)$, $({x}_1,x''_2,x_3)$$=(x''_2)$, $({x}_1,x_2,x''_3)$$=(x''_3)$. Then, for the terms in (40), (41), (42), the following estimate hold
\begin{equation}
	 |\int_{-\delta}^{\delta} \int_{\hat{x}_i}^{x'_i} \frac{\p u}{\p x_j} (x''_i)dx''_i dx'_i|\leq C,\quad 1\leq i,j\leq 3,\quad i\neq j,
\end{equation}
\begin{equation}
	 |\int_{-2\delta}^{2\delta}   \phi'(x'_i) \frac{\p u}{\p x_j} (x'_i)dx'_i|\leq C,\quad 1\leq i,j\leq 3,\quad i\neq j,
\end{equation}
\begin{equation}
	 |\int_{-2\delta}^{2\delta}   \phi (x'_i) \frac{\p u}{\p x_j} (x'_i)dx'_i|\leq C,\quad 1\leq i,j\leq 3,\quad i\neq j,
\end{equation}
\begin{equation}
	 |\int_{-\delta}^{\delta}   \frac{\p u}{\p x_j} (x'_i)dx'_i|\leq C,\quad 1\leq i,j\leq 3,\quad i\neq j,
\end{equation}
\end{theorem}

$\quad$ The proofs of Lemmas 3.1, 3.2 will be given below. Here, we admit them and give the proof of Theorem 1.2. \\

$\quad$ By inserting the estimates (43)-(46) in Lemma 3.2 into (40)-(42), we have the following. 
\begin{equation}
	(A_{11}\frac{\p u}{\p x_1})(\hat{x}_1,x_2,x_3)\leq (resp.\geq)
\end{equation}
$$
	 - (A_{12}\frac{\p u}{\p x_2}+ A_{13}\frac{\p u}{\p x_3})(\hat{x}_1,{x}_2,x_3)+(resp.-)\quad C \quad \forall (\hat{x}_1,{x}_2,x_3)\in K,
$$
\begin{equation}
	(A_{22}\frac{\p u}{\p x_2})({x}_1,\hat{x}_2,x_3)\leq (resp.\geq)
\end{equation}
$$
	 - (A_{21}\frac{\p u}{\p x_1}+ A_{23}\frac{\p u}{\p x_3})({x}_1,\hat{x}_2,x_3)+(resp.-)\quad C \quad \forall ({x}_1,\hat{x}_2,x_3)\in K,
$$
\begin{equation}
	(A_{33}\frac{\p u}{\p x_3})({x}_1,x_2,\hat{x}_3)\leq (resp.\geq)
\end{equation}
$$
	 - (A_{31}\frac{\p u}{\p x_1}+ A_{32}\frac{\p u}{\p x_2})({x}_1,{x}_2,\hat{x}_3)+(resp.-)\quad C \quad \forall (\hat{x}_1,{x}_2,\hat{x}_3)\in K. 
$$
From (47), for any $({x}_1,{x}_2,\hat{x}_3)$, $({x}_1,\hat{x}_2,x_3)$$\in K$, 
$$
	-(A_{31}\frac{\p u}{\p x_1})({x}_1,x_2,\hat{x}_3)= -(\frac{A_{31}}{A_{11}}A_{11} \frac{\p u}{\p x_1})({x}_1,x_2,\hat{x}_3)
$$
$$\leq (resp.\geq)
	 (\frac{A_{31} A_{12}}{A_{11}} \frac{\p u}{\p x_2}+\frac{A_{31}A_{13}}{A_{11}} \frac{\p u}{\p x_3})({x}_1,{x}_2,\hat{x}_3)+(resp.-)\quad C,
$$
$$
	-(A_{21}\frac{\p u}{\p x_1})({x}_1,\hat{x}_2,{x}_3)= -(\frac{A_{21}}{A_{11}}A_{11} \frac{\p u}{\p x_1})({x}_1,\hat{x}_2,{x}_3)
$$
$$\leq (resp.\geq)
	 (\frac{A_{21} A_{12}}{A_{11}} \frac{\p u}{\p x_2}+\frac{A_{21}}{A_{13}}A_{11} \frac{\p u}{\p x_3})({x}_1,\hat{x}_2,{x}_3)+(resp.-)\quad C. 
$$
Introducing the above inequalities into (48), (49) we have the following. 
\begin{equation}
	(\frac{A_{11} A_{22}- A_{21} A_{12}}{A_{11}} \frac{\p u}{\p x_2})({x}_1,\hat{x}_2,{x}_3)
\end{equation}
$$
	\leq (resp.\geq) \quad (\frac{A_{21}A_{13}-A_{11}A_{23}}{A_{11}} \frac{\p u}{\p x_3})({x}_1,\hat{x}_2,{x}_3)+(resp.-)\quad C,
$$
\begin{equation}
	(\frac{A_{11} A_{33}- A_{31} A_{13}}{A_{11}} \frac{\p u}{\p x_3})({x}_1,{x}_2,\hat{x}_3)
\end{equation}
$$
	\leq (resp.\geq) \quad (\frac{A_{31}A_{12}-A_{11}A_{32}}{A_{11}} \frac{\p u}{\p x_2})({x}_1,{x}_2,\hat{x}_3)+(resp.-)\quad C,
$$
From (50), 
$$
	(\frac{A_{31} A_{12}- A_{11} A_{32}}{A_{11}} \frac{\p u}{\p x_2})({x}_1,{x}_2,\hat{x}_3)
$$
$$
	\frac{(A_{31} A_{12}- A_{11} A_{32})   (A_{11} A_{22}- A_{21} A_{12})  }{ (A_{11} A_{22}- A_{21} A_{12})  A_{11}} \frac{\p u}{\p x_2})({x}_1,{x}_2,\hat{x}_3)
$$
$$
	\leq (resp.\geq)
$$
$$
	\frac{(A_{31} A_{12}- A_{11} A_{32})   (A_{21} A_{13}- A_{11} A_{23})  }{ (A_{11} A_{22}- A_{21} A_{12})  A_{11}} \frac{\p u}{\p x_3})({x}_1,{x}_2,\hat{x}_3)
	+(resp.-)\quad C,
$$
and by introducing the above inequalities into (51), we have the following. 
$$
	\{(A_{11} A_{33}- A_{31} A_{13})(A_{11} A_{22}- A_{21} A_{12}) - (A_{31} A_{12}- A_{11} A_{32})(A_{21} A_{13}- A_{11} A_{23}) 
	\}\times
$$
$$
	\times \frac{\p u}{\p x_3} ({x}_1,{x}_2,\hat{x}_3) \leq (resp. \geq)\quad + (resp. -)\quad C. 
$$
Therefore, from the assumptions (10), (12), we get the bound for $\frac{\p u}{\p x_3}$. A similar argument leads to the bounds for $\frac{\p u}{\p x_1}$, $\frac{\p u}{\p x_2}$, and we have proved (20).\\

{\bf Proof of Lemma 3.1.}\\
$\quad$ We only prove (40); (41), (42) will be obtained in a similar way. First of all, from (17) for any $({x}_1,{x}_2,{x}_3)$, $(\hat{x}_1,{x}_2,{x}_3) $$\in K$, since 
$$
	-(A_{11}\frac{\p u}{\p x_1} +A_{12}  \frac{\p u}{\p x_2}+A_{13} \frac{\p u}{\p x_3} ) ({x}_1,{x}_2,{x}_3) 
	+\int_{\hat{x}_1}^{x_1} (b_{12}  \frac{\p u}{\p x_2}+b_{13} \frac{\p u}{\p x_3} ) ({x}'_1,{x}_2,{x}_3) dx'_1 
$$
$$
	=-(A_{11}\frac{\p u}{\p x_1} +A_{12}  \frac{\p u}{\p x_2}+A_{13} \frac{\p u}{\p x_3} ) (\hat{x}_1,{x}_2,{x}_3) 
$$
$$
	+ \int_{\hat{x}_1}^{x_1} \{ -\frac{\p}{\p x'_1} (A_{11}\frac{\p u}{\p x'_1} +A_{12}  \frac{\p u}{\p x_2}+A_{13} \frac{\p u}{\p x_3} ) ({x}'_1,{x}_2,{x}_3) 
$$
$$
	+ (b_{12}  \frac{\p u}{\p x_2}+b_{13} \frac{\p u}{\p x_3} ) ({x}'_1,{x}_2,{x}_3) dx'_1, 
$$
by using (18), (9), (16), the following holds with $\phi$ stated in Lemma 3.1. 
$$
	|-(A_{11}\frac{\p u}{\p x_1} +A_{12}  \frac{\p u}{\p x_2}+A_{13} \frac{\p u}{\p x_3} ) ({x}_1,{x}_2,{x}_3)     
	+  \int_{\hat{x}_1}^{x_1}  (b_{12}  \frac{\p u}{\p x_2}+b_{13} \frac{\p u}{\p x_3} ) ({x}'_1,{x}_2,{x}_3) dx'_1
$$
$$
	+ (A_{11}\frac{\p u}{\p x_1} +A_{12}  \frac{\p u}{\p x_2}+A_{13} \frac{\p u}{\p x_3} ) (\hat{x}_1,{x}_2,{x}_3) |
$$
$$
	\leq \hbox{sgn}(x_1-\hat{x}_1) \int_{\hat{x}_1}^{x_1} |-\frac{\p }{\p x'_1}(A_{11}\frac{\p u}{\p x_1} +A_{12}  \frac{\p u}{\p x_2}+A_{13} \frac{\p u}{\p x_3} ) ({x}'_1,{x}_2,{x}_3)     
$$
$$
	+ (b_{11}  \frac{\p u}{\p x'_1}+ b_{12}  \frac{\p u}{\p x_2}+b_{13} \frac{\p u}{\p x_3} ) ({x}'_1,{x}_2,{x}_3) |dx'_1+C
$$
$$
	\leq \hbox{sgn}(x_1-\hat{x}_1) \int_{\hat{x}_1}^{x_1} |C_0 + \frac{\p }{\p x'_1}(A_{11}\frac{\p u}{\p x'_1} +A_{12}  \frac{\p u}{\p x_2}+A_{13} \frac{\p u}{\p x_3} ) ({x}'_1,{x}_2,{x}_3)     
$$
$$
	-(b_{11}  \frac{\p u}{\p x'_1}+ b_{12}  \frac{\p u}{\p x_2}+b_{13} \frac{\p u}{\p x_3} ) ({x}'_1,{x}_2,{x}_3) -C_0|dx'_1+C
$$
$$
	\leq \hbox{sgn}(x_1-\hat{x}_1) \int_{\hat{x}_1}^{x_1} C_0 + \frac{\p }{\p x'_1}(A_{11}\frac{\p u}{\p x'_1} +A_{12}  \frac{\p u}{\p x_2}+A_{13} \frac{\p u}{\p x_3} ) ({x}'_1,{x}_2,{x}_3)     
$$
$$
	-(b_{11}  \frac{\p u}{\p x'_1}+ b_{12}  \frac{\p u}{\p x_2}+b_{13} \frac{\p u}{\p x_3} ) ({x}'_1,{x}_2,{x}_3) dx'_1+C
$$
$$
	=\hbox{sgn}(x_1-\hat{x}_1) \times 
$$
$$
	\times  \int_{\hat{x}_1}^{x_1} \phi(x'_1) \{  C_0 + \frac{\p }{\p x'_1}(A_{11}\frac{\p u}{\p x'_1} +A_{12}  \frac{\p u}{\p x_2}+A_{13} \frac{\p u}{\p x_3} ) ({x}'_1,{x}_2,{x}_3)     
$$
$$
	-(b_{11}  \frac{\p u}{\p x'_1}+ b_{12}  \frac{\p u}{\p x_2}+b_{13} \frac{\p u}{\p x_3} ) ({x}'_1,{x}_2,{x}_3) \}dx'_1+C
$$
$$
	\leq  \int_{-2\delta}^{2\delta} \phi(x'_1) \{ \frac{\p }{\p x'_1}(A_{11}\frac{\p u}{\p x'_1} +A_{12}  \frac{\p u}{\p x_2}+A_{13} \frac{\p u}{\p x_3} ) ({x}'_1,{x}_2,{x}_3)     
$$
$$
	-(b_{11}  \frac{\p u}{\p x'_1}+ b_{12}  \frac{\p u}{\p x_2}+b_{13} \frac{\p u}{\p x_3} ) ({x}'_1,{x}_2,{x}_3) dx'_1+C
$$
$$
	\leq  \int_{-2\delta}^{2\delta} \phi'(x'_1) (A_{11}\frac{\p u}{\p x'_1} +A_{12}  \frac{\p u}{\p x_2}+A_{13} \frac{\p u}{\p x_3} ) ({x}'_1,{x}_2,{x}_3)dx'_1
$$
$$
	-  \int_{-2\delta}^{2\delta}  \phi(x'_1) (b_{12}  \frac{\p u}{\p x_2}+b_{13} \frac{\p u}{\p x_3} ) ({x}'_1,{x}_2,{x}_3) \}dx'_1+C
$$
$$
	=- \int_{-2\delta}^{2\delta} \phi'(x'_1) (A_{12}  \frac{\p u}{\p x_2}+A_{13} \frac{\p u}{\p x_3} ) ({x}'_1,{x}_2,{x}_3)dx'_1
$$
$$
	-  \int_{-2\delta}^{2\delta}  \phi(x'_1) (b_{12}  \frac{\p u}{\p x_2}+b_{13} \frac{\p u}{\p x_3} ) ({x}'_1,{x}_2,{x}_3) \}dx'_1+C
$$
Hence, for any $({x}_1,{x}_2,{x}_3) $, $(\hat{x}_1,{x}_2,{x}_3) $$\in K$, we get the following inequalities. 
\begin{equation}
	(A_{11}\frac{\p u}{\p x_1} +A_{12}  \frac{\p u}{\p x_2}+A_{13} \frac{\p u}{\p x_3})(x_1,x_2,x_3)
\end{equation}
$$
	\leq (resp.\geq) \quad (A_{11}\frac{\p u}{\p x_1} +A_{12}  \frac{\p u}{\p x_2}+A_{13} \frac{\p u}{\p x_3})(\hat{x}_1,x_2,x_3)
$$
$$
	+ \int_{\hat{x}_1}^{x_1}  (b_{12}  \frac{\p u}{\p x_2}+b_{13} \frac{\p u}{\p x_3} ) ({x}'_1,{x}_2,{x}_3) dx'_1
$$
$$
	- (resp.+)\quad \int_{-2\delta}^{2\delta} \phi'(x'_1) (A_{12}  \frac{\p u}{\p x_2}+A_{13} \frac{\p u}{\p x_3} ) ({x}'_1,{x}_2,{x}_3)dx'_1
$$
$$
	-(resp.+)\quad  \int_{-2\delta}^{2\delta} \phi (x'_1) (b_{12}  \frac{\p u}{\p x_2}+b_{13} \frac{\p u}{\p x_3} ) ({x}'_1,{x}_2,{x}_3) dx'_1+
	(resp.-)\quad C
$$
$\quad$ Next, we integrate the both hands sides of the above inequalities with respect to $x_1$ on $[-\delta,\delta]$, then devide the result by $2\delta$ and we have the following. 
\begin{equation}
	\frac{1}{2\delta} \int_{-\delta}^{\delta} (A_{12}  \frac{\p u}{\p x_2}+A_{13} \frac{\p u}{\p x_3})(x'_1,x_2,x_3)dx'_1
\end{equation}
$$
	\leq (resp.\geq)\quad (A_{11}\frac{\p u}{\p x_1} +A_{12}  \frac{\p u}{\p x_2}+A_{13} \frac{\p u}{\p x_3})(\hat{x}_1,x_2,x_3)
$$
$$
	+\frac{1}{2\delta} \int_{-\delta}^{\delta} \int_{\hat{x}_1}^{x_1}  (b_{12}  \frac{\p u}{\p x_2}+b_{13} \frac{\p u}{\p x_3} ) ({x}''_1,{x}_2,{x}_3) \}dx''_1dx'_1
$$
$$
	-(resp.+)\quad  \int_{-2\delta}^{2\delta}  \phi' (x'_1) (A_{12}  \frac{\p u}{\p x_2}+A_{13} \frac{\p u}{\p x_3} ) ({x}'_1,{x}_2,{x}_3) \}dx'_1
$$
$$
	-(resp.+)\quad  \int_{-2\delta}^{2\delta}  \phi (x'_1) (b_{12}  \frac{\p u}{\p x_2}+b_{13} \frac{\p u}{\p x_3} ) ({x}'_1,{x}_2,{x}_3) \}dx'_1+
	(resp.-)\quad C
$$
The above inequality leads (40).\\

{\bf Proof of Lemma 3.2.}\\
$\quad$ We show the estimates (43)-(46) in the following steps 1-4.\\

{\bf Step 1.} (Estimate (43).) We consider the particular case when $i=1$, $j=2$; the other cases are obtained in a similar way in view of the symmetry of the conditions on $(A_{ij})$ and $(b_{ij})$. 
\begin{equation}
	|\int_{-\delta}^{\delta} \int_{\hat{x}_1}^{x'_1}  \frac{\p u}{\p x_2} ({x}''_1,\hat{x}_2,{x}_3) dx''_1dx'_1|\leq C,\quad \forall \hat{x}_2,,\forall {x}_3\in [-\delta,\delta]. 
\end{equation}
First, we integrate both hands sides of the inequalities (41) with respect to $x_1$ on $[\hat{x}_1, x'_1]$, and then with respect to $x'_1$ on $[-\delta,\delta]$. 
Remark that since $u$ is bounded, the integrals of   $\frac{\p u}{\p x_3}$ $({x}_1,\hat{x}_2,{x}_3)$ with respect to $x_1$ are estimated by constants. Moreover, remark that by using (42) and the boundedness of $u$, the integrals of $\frac{\p u}{\p x_3} ({x}_1,{x}_2,{x}_3)$ with respect to $x_2$ and $x_1$ are estimated by constants. Thus, we get the following inequality from (41) and (42). 
$$
	\int_{-\delta}^{\delta} \int_{\hat{x}_1}^{x'_1}  A_{22} \frac{\p u}{\p x_2} ({x}''_1,\hat{x}_2,{x}_3) dx''_1dx'_1
$$
$$(resp.\geq)\quad \int_{-\delta}^{\delta} \int_{\hat{x}_1}^{x'_1}  -A_{23} \frac{\p u}{\p x_3} ({x}''_1,\hat{x}_2,{x}_3) dx''_1dx'_1
	+(resp.-)\quad C.
$$
We denote
$$
	B(x_3)=\int_{-\delta}^{\delta}\int_{\hat{x}_1}^{x'_1}  A_{22} \frac{\p u}{\p x_2} ({x}''_1,\hat{x}_2,{x}_3) dx''_1dx'_1,
$$
where $\hat{x}_2$$\in [-\delta,\delta]$ is arbitrarily fixed. Then, since 
$$
	 -A_{23} \frac{\p u}{\p x_3} ({x}''_1,\hat{x}_2,{x}_3)= -\frac{A_{23}}{A_{33}}A_{33} \frac{\p u}{\p x_3} ({x}''_1,\hat{x}_2,{x}_3), 
$$
by inserting (42) into the above inequalities and by using the boundedness of $u$, we deduce 
\begin{equation}
	A_{22}B(x_3)\leq (resp.\geq)\quad \frac{A_{23}A_{32}}{A_{33}} B(x_3) - \frac{1}{2\delta} \frac{A_{23}A_{32}}{A_{33}} \int_{-\delta}^{\delta} B(x'_3) dx'_3
\end{equation}
$$
	+ \frac{1}{2\delta} \frac{A_{23}b_{32}}{A_{33}} \int_{-\delta}^{\delta} \int_{\hat{x}_3}^{x'_3} B(x''_3) dx''_3dx'_3
$$
$$
	-(resp.+) |\frac{A_{23}}{A_{33}} | A_{32} \int_{-2\delta}^{2\delta} \phi'(x'_3)B(x'_3) dx'_3
$$
$$
	-(resp.+)\quad  |\frac{A_{23}}{A_{33}} |b_{32} \int_{-2\delta}^{2\delta} \phi'(x'_3)B(x'_3) dx'_3+(resp.-)\quad C.
$$
We multiply (55) by $\phi'(x_3)$ and integrate the result with respect to $x_3$ on $[-2\delta,2\delta]$. Then, from the assumption on $\phi$ in (39), we have 
$$
	\frac{A_{22}A_{33}-A_{23}A_{32}}{A_{33}}\int_{-2\delta}^{2\delta} \phi'(x'_3) B(x'_3) dx'_3
$$
$$
	\leq (resp.\geq)\quad -(resp.+) 2 |\frac{A_{23}}{A_{33}} |A_{32} \int_{-2\delta}^{2\delta} \phi'(x'_3)B(x'_3) dx'_3
$$
$$
	 -(resp.+) 2 |\frac{A_{23}}{A_{33}} |b_{32} \int_{-2\delta}^{2\delta} \phi(x'_3)B(x'_3) dx'_3+(resp.-)\quad C.
$$
Here, we shall denote 
\begin{equation}
	E_1=(\frac{A_{22}A_{33}-A_{23}A_{32}}{A_{33}} -2|\frac{A_{23}}{A_{33}}|A_{32})^{-1} |\frac{A_{23}}{A_{33}}|b_{32}
\end{equation}

\begin{equation}
	E_2=(\frac{A_{22}A_{33}A_{23}A_{32}}{A_{33}} +2|\frac{A_{23}}{A_{33}}|A_{32})^{-1} |\frac{A_{23}}{A_{33}}|b_{32}
\end{equation}
From the condition (15), we have the following two cases. \\

Case 1. The following inequalities hold.
$$
	\frac{A_{22}A_{33}-A_{23}A_{32}}{A_{33}} +2|\frac{A_{23}}{A_{33}}|A_{32}>0, 
$$
\begin{equation}
\quad
\end{equation}
$$
	\frac{A_{22}A_{33}-A_{23}A_{32}}{A_{33}} -2|\frac{A_{23}}{A_{33}}|A_{32}>0. 
$$
Case 2. The following inequalities hold. 
$$
	\frac{A_{22}A_{33}-A_{23}A_{32}}{A_{33}} +2|\frac{A_{23}}{A_{33}}|A_{32}<0, 
$$
\begin{equation}
\quad
\end{equation}
$$
	\frac{A_{22}A_{33}-A_{23}A_{32}}{A_{33}} -2|\frac{A_{23}}{A_{33}}|A_{32}<0. 
$$
So, in Case 1 ((58)), 
$$
	\int_{-2\delta}^{2\delta} \phi'(x'_3) B(x'_3)dx'_3 \leq -2E_2\int_{-2\delta}^{2\delta} \phi(x'_3) B(x'_3)dx'_3, 
$$
$$
	\int_{-2\delta}^{2\delta} \phi'(x'_3) B(x'_3)dx'_3 \geq 2E_1 \int_{-2\delta}^{2\delta} \phi(x'_3) B(x'_3)dx'_3, 
$$
and in Case 2 ((59)), 
$$
	\int_{-2\delta}^{2\delta} \phi'(x'_3) B(x'_3)dx'_3 \geq -2E_2\int_{-2\delta}^{2\delta} \phi(x'_3) B(x'_3)dx'_3, 
$$
$$
	\int_{-2\delta}^{2\delta} \phi'(x'_3) B(x'_3)dx'_3 \leq 2E_1 \int_{-2\delta}^{2\delta} \phi(x'_3) B(x'_3)dx'_3. 
$$
By inserting these inequalities into (55), we get the following
\begin{equation}
	\frac{A_{22}A_{33}A_{23}A_{32}}{A_{33}}B(x_3)\leq (resp.\geq)\quad -\frac{1}{2\delta} \frac{A_{23}A_{32}}{A_{33}}\int_{-\delta}^{\delta} B(x'_3)dx'_3
\end{equation}
$$
	+\frac{1}{2\delta} \frac{A_{23}b_{32}}{A_{33}}\int_{-\delta}^{\delta} \int_{\hat{x}_3}^{x'_3} B(x''_3)dx''_3dx'_3
$$
$$
	-(resp.+)\quad \frac{A_{22}A_{33}-A_{23}A_{32}}{A_{33}} E_i \int_{-2\delta}^{2\delta} \phi(x'_3)B(x'_3)dx'_3+(resp.-)\quad C,
$$
where $i=1$ in Case 1 and $A_{32}\geq 0$, or in Case 2 and $A_{32}\leq 0$; $i=2$ in Case 1 and $A_{32}\leq 0$, or in Case 2 and $A_{32}\geq 0$. \\
$\quad$ Next, we investigate both hand sides of (60) with respect to $x_3$ on $[-\delta,\delta]$, and devide both hands sides of the result by $A_{22}$.   
\begin{equation}
	\int_{-\delta}^{\delta} B(x'_3)dx'_3 \leq (resp.\geq)\quad \frac{A_{23}b_{32}}{A_{22}A_{33}}\int_{-\delta}^{\delta} \int_{\hat{x}_3}^{x'_3} B(x''_3)dx''_3dx'_3
\end{equation}
$$
	-(resp.+)\quad \frac{A_{22}A_{33}-A_{23}A_{32}}{A_{33}} \frac{2\delta E_i}{|A_{22}|} \int_{-2\delta}^{2\delta} \phi(x'_3)B(x'_3)dx'_3+(resp.-)\quad C,
$$
where the indices $i=1,2$ are similar to (60). By inserting the above inequalities into (60), then deviding both hands sides of the result by $\alpha=$$\frac{A_{22}A_{33}-A_{23}A_{32}}{A_{33}} $, we get the following. 
\begin{equation}
	B(x_3) \leq (resp.\geq)\quad\frac{1}{2\delta} \frac{A_{23}b_{32}}{A_{22}A_{33}}\int_{-\delta}^{\delta} \int_{\hat{x}_3}^{x'_3} B(x''_3)dx''_3dx'_3
\end{equation}
$$
	-(resp.+)\quad (\hbox{sgn}\alpha) E_i (1+|\frac{A_{23}A_{32}}{A_{22}A_{33}} |)\int_{-2\delta}^{2\delta} \phi(x'_3)B(x'_3)dx'_3+(resp.-)\quad C,
$$
where the indices $i=1,2$ are similar to (60). \\
$\quad$ We integrate the both hand sides of (62) with respect to $x_3$ on $[\hat{x}_3,x'_3]$, then with respect to $x'_3$ on $[-\delta,\delta]$, which leads the following. 
$$
	\int_{-\delta}^{\delta} \int_{\hat{x}_3}^{x'_3}B(x''_3)dx''_3dx'_3 \leq (resp.\geq)\quad
	-\hat{x}_3\frac{A_{23}b_{32}}{A_{22}A_{33}}\int_{-\delta}^{\delta} \int_{\hat{x}_3}^{x'_3} B(x''_3)dx''_3dx'_3
$$
$$
	-(resp.+)\quad (-2\delta \hat{x}_3) (\hbox{sgn}\alpha) E_i (1+\frac{A_{23}A_{32}}{A_{22}A_{33}}|) \int_{-2\delta}^{2\delta} \phi(x'_3)B(x'_3)dx'_3+(resp.-)\quad C,
$$
where the indices $i=1,2$ are similar to (60). From (15), 
$$
	\frac{1}{A_{22}A_{33}}\leq \max\{ \frac{1}{A_{22}A_{33}-A_{23}A_{32}+2|A_{23}A_{32}|}, 
	\frac{1}{A_{22}A_{33}-A_{23}A_{32}-2|A_{23}A_{32}|}
	\}, 
$$
and since $|\hat{x}_3|<\delta$, we have from (19)  
$$
	|-\hat{x}_3\frac{A_{23}b_{32}}{A_{22}A_{33}}|\leq \frac{1}{2}. 
$$
Thus, for each cases of $i=1,2$, there exist constants $O_1(\delta)$$=O(\delta)$, $O_2(\delta)$$=O(\delta)$ respectively, such that 
$$
	\frac{1}{2\delta} \frac{A_{23}b_{32}}{A_{22}A_{33}}\int_{-\delta}^{\delta} \int_{\hat{x}_3}^{x'_3}B(x''_3)dx''_3dx'_3
$$
$$
	\leq(resp.\geq)\quad -(resp.+)\quad O_i(\delta) \int_{-2\delta}^{2\delta} \phi(x'_3)B(x'_3)dx'_3 +(resp.-)\quad C,
$$
where the indices $i=1,2$ are similar to (60). By inserting the above estimate into (62), we get 
\begin{equation}
	B(x_3) \leq(resp.\geq)
\end{equation}
$$
	-(resp.+)\quad \{O_i(\delta)+  (\hbox{sgn}\alpha) E_i (1+|\frac{A_{23}A_{32}}{A_{22}A_{33}} | \}\int_{-2\delta}^{2\delta} \phi(x'_3)B(x'_3)dx'_3 +(resp.-)\quad C,
$$
where the indices $i=1,2$ are similar to (60). We multiply (63) by $\phi(x_3)>0$ and integrate both hand sides of the result with respect to $x_3$ on $[-2\delta,2\delta]$. Then, by remarking that 
$$
	2\delta\leq \int_{-2\delta}^{2\delta} \phi(x'_3)B(x'_3)dx'_3 \leq 4\delta,
$$
also by remarking that from (19), 
$$
	|4\delta E_i (1+|\frac{A_{23}A_{32}}{A_{22}A_{33}}|) |<\frac{1}{2} \quad i=1,2, 
$$
and by noticing that $O_i(\delta)$$=O(\delta)$ for $i=1,2$, we obtain 
$$
	|\int_{-2\delta}^{2\delta} \phi(x'_3)B(x'_3)dx'_3|<C.
$$
By inserting the last estimate into (63), we obtain the estimate. \\

{\bf Step 2.} (Estimate (44).)\\
$\quad$ We consider the particular case when $i=1$, $j=2$, by separating it into the following two inequalities; the other cases are obtained in a similar way 
in view of the symmetry of the conditions on $(A_{ij})_{1\leq i,j\leq 3}$ and  $(b_{ij})_{1\leq i,j\leq 3}$. 
\begin{equation}
	|\int_{-2\delta}^{0} \phi'(x'_1,\hat{x}_2,x_3)  \frac{\p u}{\p x_2} (x'_1,\hat{x}_2,x_3) dx'_1|<C, \quad \forall \hat{x}_2,\quad x_3\in [-\delta,\delta], 
\end{equation}
\begin{equation}
	|\int_{0}^{2\delta} \phi'(x'_1,\hat{x}_2,x_3)  \frac{\p u}{\p x_2} (x'_1,\hat{x}_2,x_3) dx'_1|<C, \quad \forall \hat{x}_2,\quad x_3\in [-\delta,\delta], 
\end{equation}
It is enough to show (64), because (65) can be proved in the same way. Now, we set 
$$
	C(x_3)= \int_{-2\delta}^{0} \phi'(x'_1,\hat{x}_2,x_3)  \frac{\p u}{\p x_2} (x'_1,\hat{x}_2,x_3) dx'_1, 
$$
where $\hat{x}_2$$\in [-\delta,\delta]$ is arbitrarily fixed. By using the estimate (43) in (41), in the same way as in Step 1, we obtain 
\begin{equation}
	A_{22}C(x_3) \leq (resp.\geq)\quad \frac{A_{23}A_{32}}{A_{33}}C(x_3) - \frac{1}{2\delta} \frac{A_{23}A_{32}}{A_{33}} \int_{-\delta}^{\delta}C(x'_3) dx'_3
\end{equation}
$$
	-(resp.+)\quad |\frac{A_{23}}{A_{33}}| A_{32} \int_{-2\delta}^{2\delta} \phi'(x'_3)C(x'_3)dx'_3
$$
$$
	-(resp.+)\quad |\frac{A_{23}}{A_{33}}| b_{32} \int_{-2\delta}^{2\delta} \phi(x'_3)C(x'_3)dx'_3+(resp.-) \quad C. 
$$
We multiply both hand sides of (66) by $\phi'(x_3)$, then integrate the result with respect to $x_3$ on $ [-2\delta,2\delta]$. From the assumption on $\phi$ in (39), we get 
$$
	\frac{A_{22}A_{33}-A_{23}A_{32}}{A_{33}}\int_{-2\delta}^{2\delta} \phi'(x'_3)C(x'_3)dx'_3
$$
$$
	\leq (resp.\geq)\quad -(resp.+)\quad 2 |\frac{A_{23}}{A_{33}}| A_{32} \int_{-2\delta}^{2\delta} \phi'(x'_3)C(x'_3)dx'_3
$$
$$
	-(resp.+)\quad 2 |\frac{A_{23}}{A_{33}}| b_{32} \int_{-2\delta}^{2\delta} \phi(x'_3)C(x'_3)dx'_3+(resp.-) \quad C. 
$$
From the condition (15), we have the following two cases. \\
Case 1.\\
The following inequalities hold. 
$$
	\frac{A_{22}A_{33}-A_{23}A_{32}}{A_{33}}+ 2|\frac{A_{23}}{A_{33}}| A_{32}>0, 
$$
$$
	\frac{A_{22}A_{33}-A_{23}A_{32}}{A_{33}}- 2|\frac{A_{23}}{A_{33}}| A_{32}>0. 
$$
Case 2.\\
The following inequalities hold. 
$$
	\frac{A_{22}A_{33}-A_{23}A_{32}}{A_{33}}+ 2|\frac{A_{23}}{A_{33}}| A_{32}<0, 
$$
$$
	\frac{A_{22}A_{33}-A_{23}A_{32}}{A_{33}}- 2|\frac{A_{23}}{A_{33}}| A_{32}<0. 
$$
Thus, denoting by 
\begin{equation}
	E_1= (\frac{A_{22}A_{33}-A_{23}A_{32}}{A_{33}}- 2|\frac{A_{23}}{A_{33}}| A_{32})^{-1}|\frac{A_{23}}{A_{33}}| b_{32},
\end{equation}
\begin{equation}
	E_2= (\frac{A_{22}A_{33}-A_{23}A_{32}}{A_{33}}+ 2|\frac{A_{23}}{A_{33}}| A_{32})^{-1}|\frac{A_{23}}{A_{33}}| b_{32},
\end{equation}
the same argument as in Step 1 to deduce (61) leads to the following inequalities. 
\begin{equation}
	\frac{A_{22}A_{33}-A_{23}A_{32}}{A_{33}}C(x_3)\leq (resp.\geq)\quad -\frac{1}{2\delta} \frac{A_{23}A_{32}}{A_{33}}\int_{-\delta}^{\delta} C(x'_3)dx'_3
\end{equation}
$$
	-(resp.+)\quad \frac{A_{22}A_{33}-A_{23}A_{32}}{A_{33}} E_i \int_{-2\delta}^{2\delta} \phi(x'_3)C(x'_3)dx'_3 +(resp.-) \quad C,
$$
where $i=1$ in Case 1 and $A_{32}\geq 0$, or in Case 2 and $A_{32}\leq 0$; $i=2$ in Case 1 and $A_{32}\leq 0$, or in Case 3 and $A_{32}\geq 0$. \\
$\quad$ Next, by integrating both hand sides of (69) with respect to $x_3$ on $[-\delta,\delta]$, then by deviding  the result by $A_{22}$, we get the following.  
$$
	\int_{-\delta}^{\delta} C(x'_3)dx'_3 \leq (resp.\geq) 
$$
$$
	-(resp.+)\quad  2\delta  E_i  \frac{A_{22}A_{33}-A_{23}A_{32}}{|A_{22}|A_{33}} \int_{-2\delta}^{2\delta} \phi(x'_3)C(x'_3)dx'_3 +(resp.-) \quad C,
$$
where the indices $i=1,2$ are similar to (69). By inserting this inequality into (69) and devide both hand sides of the result by $\alpha=$$\frac{A_{22}A_{33}-A_{23}A_{32}}{|A_{22}|A_{33}}$, we get the following. 
\begin{equation}
	C(x_3)\leq (resp.\geq) 
\end{equation}
$$
	-(resp.+)\quad  (\hbox{sgn}\alpha) E_i  (1+|\frac{A_{23}A_{32}}{A_{22}A_{33}}|) \int_{-2\delta}^{2\delta} \phi(x'_3)C(x'_3)dx'_3 +(resp.-) \quad C. 
$$
By remarking that from (19), 
$$
	|4\delta  E_i(1+|\frac{A_{23}A_{32}}{A_{22}A_{33}}|)|<\frac{1}{2} \quad i=1,2,
$$
and by using the same argument as in Step 1, we have 
$$
	|\int_{-2\delta}^{2\delta} \phi(x'_3)C(x'_3)dx'_3|<C. 
$$
By inserting the last estimate into (70), we obtain the estimate (64); (65) can be obtained by the same way. \\

{\bf Step 3.} (Estimate (45).)\\
$\quad$ We consider the particular case when $i=1$, $j=2$; the other cases are obtained in a similar way in view of the symmetry of the coefficients $(A_{ij})$ and $(b_{ij})$. 
\begin{equation}
	|\int_{-2\delta}^{2\delta} \phi(x'_1) \frac{\p u}{\p x_2} (x'_1) dx'_1|<C. 
\end{equation}
We set
$$
	D(x_3)= \int_{-2\delta}^{2\delta} \phi(x'_1,\hat{x}_2,x_3)  \frac{\p u}{\p x_2} (x'_1,\hat{x}_2,x_3) dx'_1, 
$$
where $\hat{x}_2$$\in [-\delta,\delta]$ is arbitrarily fixed. By inserting the estimates (43), (44) into (41), and by using the same argument as in Steps 1,2, we get the following 
\begin{equation}
	A_{22}D(x_3) \leq (resp.\geq)\quad \frac{A_{23}A_{32}}{A_{33}}D(x_3) -\frac{1}{2\delta} \frac{A_{23}A_{32}}{A_{33}} \int_{-\delta}^{\delta}D(x'_3)dx'_3
\end{equation}
$$
	-(resp.+)\quad |\frac{A_{23}}{A_{33}}|b_{32} \int_{-2\delta}^{2\delta} \phi(x'_3)D(x'_3)dx'_3+(resp.-)\quad C. 
$$
By integrating both hand sides of (72) with respect to $x_3$ on $[-\delta,\delta]$, then by deviding the result by $A_{22}$, we have the following. 
\begin{equation}
	\int_{-\delta}^{\delta}D(x'_3)dx'_3\leq (resp.\geq)
\end{equation}
$$
	-(resp.+) \quad 2\delta |\frac{A_{23}A_{32}}{A_{33}}|b_{32}  \int_{-2\delta}^{2\delta} \phi(x'_3)D(x'_3)dx'_3+(resp.-)\quad C. 
$$
By inserting (73) into (72), we obtain 
$$
	 \frac{A_{22}A_{33}-A_{23}A_{32}}{A_{33}}D(x_3)  \leq (resp.\geq) 
$$
$$
	-(resp.+) \quad |\frac{A_{23}}{A_{33}}|b_{32}(1+ |\frac{A_{23}A_{32}}{A_{22}A_{33}}| )\int_{-2\delta}^{2\delta} \phi(x'_3)D(x'_3)dx'_3+(resp.-)\quad C, 
$$
where $C>0$ is a constant. We multiply both hand sides of the above inequalities by $\phi(x_3)$, then integrate the result with respect to $x_3$ on $[-2\delta,2\delta]$. Since from the assumption (15), 
$$
	 \frac{1}{|A_{22}A_{33}-A_{23}A_{32}|}\leq 
$$
$$
	\max\{\frac{1}{A_{22}A_{33}-A_{23}A_{32}+ 2|A_{23}A_{32}|}, 
	\frac{1}{A_{22}A_{33}-A_{23}A_{32}- 2|A_{23}A_{32}|}
	\},
$$
remarking that from (19), 
$$
	4\delta |\frac{A_{22}A_{33}-A_{23}A_{32}}{A_{33}}|^{-1} |\frac{A_{23}}{A_{33}}||b_{32}| (1+ |\frac{A_{23}A_{32}}{A_{22}A_{33}}| )<\frac{1}{2}, 
$$
we get 
$$
	|\int_{-2\delta}^{2\delta} \phi(x'_3)D(x'_3)dx'_3|<C.
$$
By inserting the above estimates into (72), (73), from (11), we obtain (71). \\

{\bf Step 4.} (Estimate (46).)\\
$\quad$ We consider the particular case when $i=1$, $j=2$; the other cases are obtained in a similar way, in view of the symmetry of the coefficients $(A_{ij})$, 
$(b_{ij})$. 
\begin{equation}
	|\int_{-\delta}^{\delta} \frac{\p u}{\p x_2}(x'_1) dx'_1|<C.
\end{equation}
We set 
$$
	E(x_3)= \int_{-\delta}^{\delta} \frac{\p u}{\p x_2}(x'_1,\hat{x}_2,x_3) dx'_1,
$$
where $\hat{x}_2$ is arbitrarily fixed. Then, by inserting the estimates (43)-(45) into (41), and by taking the same arguments as in Steps 1-3, we get 
$$
	A_{22}E(x_3)\leq (resp.\geq) \quad \frac{A_{23}A_{32}}{A_{33}}E(x_3) 
	-\frac{1}{2\delta} \frac{A_{23}A_{32}}{A_{33}}  \int_{-\delta}^{\delta} E(x'_3)dx'_3 +(resp.-)\quad C, 
$$
where $C>0$ is a constant. By integrating the both hand sides of the above inequalities, we have the estimate 
$$
	| \int_{-\delta}^{\delta} E(x'_3)dx'_3  |\leq C, 
$$
and by inserting this into the above inequality, we obtain (74). \\
$\quad $ From Steps 1-4, we have proved Lemma 3.2. \\

\section{Proofs of Theorems 1.3, 1.4 and Example 1.7.}

$\quad$ We begin with the following Lemma.\\

\leftline{\bf Lemma 4.1.}
\begin{theorem}
 For the function $u$ in Theorem 1.3, the following estimates hold. 
\begin{equation}
	\sup_{x\in \Omega} | (\sum_{j=1}^{N} A_{ij} \frac{\p u}{\p x_i})(x)| \leq C\quad 1\leq i\leq N. 
\end{equation}
\end{theorem}

$\quad$ From Lemma 4.1, (21), (22), we obtain the gradient estimate (24) in Theorem 1.3. \\

Now, we prove Lemma 4.1.\\
{\bf Proof of Lemma 4.1.}\\
$\quad$ We only show the estimate for $i=1$ in (75); and the others are obtained by the same way. Let us denote $(x_1,...,x_N)$$=(x_1)$,    $(x'_1,...,x_N)$$=(x'_1)$, $(\hat{x}_1,...,x_N)$$=(\hat{x}_1)$ for the convenience. For any $(x_1)$, $(\hat{x}_1)\in \Omega$,  since 
$$
	-\sum_{j=1}^{N}( A_{1j} \frac{\p u}{\p x_j})(x_1)
$$
$$
	= -\sum_{j=1}^{N} (A_{1j} \frac{\p u}{\p x_j})(\hat{x}_1) + \int_{\hat{x}_1}^{x_1} -\frac{\p }{\p x'_1}\sum_{j=1}^{N} (A_{1j} \frac{\p u}{\p x_j})({x}'_1)dx'_1, 
$$
the following holds from (23). 
$$
	| -\sum_{j=1}^{N} A_{1j} \frac{\p u}{\p x_j})(x_1) +\sum_{j=1}^{N} A_{1j} \frac{\p u}{\p x_j})(\hat{x}_1)  |
$$
$$
	\leq \hbox{sgn}(x_1-\hat{x}_1)\int_{\hat{x}_1}^{x_1} | -\frac{\p }{\p x'_1}\sum_{j=1}^{N} (A_{1j} \frac{\p u}{\p x_j})({x}'_1)|dx'_1 
$$
$$
	= \hbox{sgn}(x_1-\hat{x}_1)\int_{\hat{x}_1}^{x_1}  | C_0
+\frac{\p }{\p x'_1}\sum_{j=1}^{N} (A_{1j} \frac{\p u}{\p x_j})({x}'_1)-C_0|dx'_1 
$$
$$
	\leq \hbox{sgn}(x_1-\hat{x}_1)\int_{\hat{x}_1}^{x_1} C_0+\frac{\p }{\p x'_1}\sum_{j=1}^{N} (A_{1j} \frac{\p u}{\p x_j})({x}'_1) dx'_1 +C, 
$$
$$
	\leq \int_y^z C_0+\frac{\p }{\p x'_1}\sum_{j=1}^{N} (A_{1j} \frac{\p u}{\p x_j})({x}'_1) dx'_1 + C, 
$$
where $(y,x_2,...,x_N)$$\in \p\Omega$, $(z,x_2,...,x_N)$$\in \p\Omega$, ($y\leq z$) are the intersections of $\p\Omega$ and the straight line connecting $(x_1)$ with $(\hat{x}_1)$. Hence, from the assumption that $\hbox{supp}u$$\subset \subset \Omega$, we get 
$$
	| -\sum_{j=1}^{N} A_{1j}\frac{\p u}{\p x_j}(x_1) +  \sum_{j=1}^{N} A_{1j}\frac{\p u}{\p x_j}(\hat{x}_1) |\leq C. 
$$
And, by letting $\hat{x}_1$ be on the boundary, we have proved our purpose. \\

$\quad$ For Theorem 1.4, the same lemma as above holds. \\

\leftline{\bf Lemma 4.2.}
\begin{theorem}
For the function $u$ in Theorem 1.4, the following estimate holds. 
\begin{equation}
	\sup_{x\in \Omega}|  (A_{ij}\frac{\p u}{\p x_i}(x) |\leq C \quad 1\leq i\leq N. 
\end{equation}
\end{theorem}

$\quad$ It is not difficult to prove Lemma 4.2, by modifying the proof of Lemma 4.1. Moreover, it is clear that Lemma 4.2 leads Theorem 1.4, and we do not repeat the argument. \\

{\bf Proof of Example 1.7.}\\
$\quad$ The existence and the unqueness of the solution $u_{\l}$ of (29) is established by the viscosity solutions theory. (We refer the viscosity solutions theory to Crandall-Lions [3], Crandall-Ishii-Lions [4].) Thus, by the comparison result, we have 
$$
	-\frac{\p^2 u_{\l}}{\p x_i^2}\leq \hbox{Const.} \quad x\in\Omega, \quad 1\leq i\leq N, 
$$
$$
	\l u_{\l} (x)\leq \hbox{Const.}
$$
Therefore, we can apply Theorems 1.1, 1.2, and 1.4 to obtain the result.\\

\leftline{\bf Remarks 4.3.} The regularity result in Example 1.7 can be generalized to a class of some controlled stochastic systems which were treated by Krylov [6], Lions [7]. For the special case of (29), the result in fact holds for any dimensions, if we follow our proof for Theorems 1.1, 1.2, and 1.4. \\

\end{document}